\documentclass{amsart}

\usepackage{amsfonts}
\usepackage{amsmath}
\usepackage{amsthm}
\usepackage{amssymb}
\usepackage{ifsym}
\usepackage{dirtytalk}

\theoremstyle{definition}

\theoremstyle{remark}

\numberwithin{equation}{section}

\begin{document}

\title[]{Decomposing the real line into everywhere isomorphic suborders}
\author{Garrett Ervin}

\begin{abstract}
We show that if $\mathbb{R} = A \cup B$ is a partition of $\mathbb{R}$ into two suborders $A$ and $B$, then there is an open interval $I$ such that $A \cap I$ is not order-isomorphic to $B \cap I$. The proof depends on the completeness of $\mathbb{R}$, and we show in contrast that there is a partition of the irrationals $\mathbb{R} \setminus \mathbb{Q} = A \cup B$ such that $A \cap I$ is isomorphic to $B \cap I$ for every open interval $I$. We do not know if there is a partition of $\mathbb{R}$ into three suborders that are isomorphic in every open interval. 
\end{abstract}

\maketitle

\section{Introduction}

Suppose that $A$ and $B$ are disjoint subsets of $\mathbb{R}$, viewed as suborders of $(\mathbb{R}, <)$. We say that $A$ and $B$ are \emph{everywhere isomorphic} if for every open interval $I = (a, b)$ we have that $A \cap I$ is order-isomorphic to $B \cap I$. We allow $a = -\infty$ and $b = \infty$ in this definition, so that if $A$ and $B$ are everywhere isomorphic, they are in particular isomorphic as linear orders.

If $A$ and $B$ are disjoint countable subsets of $\mathbb{R}$ that are each dense in $\mathbb{R}$, then $A$ and $B$ are everywhere isomorphic, since on every open interval $I$ we have that $A \cap I \cong B \cap I \cong \mathbb{Q}$. This follows from Cantor's theorem that every countable dense linear order without endpoints is order-isomorphic to the rationals \cite[pg. 122]{Cantor}. For instance, $\mathbb{Q}$ and its shift $\mathbb{Q} + \sqrt{2} = \{q + \sqrt{2}: q \in \mathbb{Q}\}$ are everywhere isomorphic. 

There are also uncountable examples of everywhere isomorphic sets. In fact, for any infinite cardinal $\kappa \leq 2^{\aleph_0}$ one can find disjoint subsets $A$ and $B$ of $\mathbb{R}$ of size $\kappa$ that are everywhere isomorphic. We describe a method for constructing such sets below, and show moreover that they can be homogeneous as linear orders. 

Given that there are examples of everywhere isomorphic sets $A$ and $B$ with $|A| = |B| = 2^{\aleph_0}$, it is natural to ask whether $\mathbb{R}$ can be partitioned into everywhere isomorphic ``half-sets," that is, whether there are two everywhere isomorphic sets whose union is $\mathbb{R}$. Our main result is that no such partition exists. 

\theoremstyle{definition}
\newtheorem{thm1}{Theorem}
\begin{thm1}\label{thm1} \,\ 
If $\mathbb{R} = A \cup B$ is a partition of $\mathbb{R}$, then there is an open interval $I$ such that $A \cap I$ is not order-isomorphic to $B \cap I$. 
\end{thm1}

There are many partitions $\mathbb{R} = A \cup B$ for which $A$ and $B$ are isomorphic globally. For example, let $A$ be the union of all half-open intervals $[2n, 2n+1)$ with $n \in \mathbb{Z}$, and let $B$ be the union of the intervals $[2n+1, 2n+2)$. Then $x \mapsto x+1$ defines an order-isomorphism from $A$ to $B$. Of course, this $A$ and $B$ are not everywhere isomorphic, since for any $n \in \mathbb{Z}$ and any open interval $I \subseteq [n, n+1)$, one of the two sets $A \cap I$ and $B \cap I$ is empty. A variation of this example gives a partition of $\mathbb{R}$ into two globally isomorphic sets that are each dense in $\mathbb{R}$, and even have the same cardinality in every open interval. But Theorem 1 shows that there is no construction that gets isomorphism in every open interval. 

The proof of Theorem \ref{thm1} uses the completeness of $(\mathbb{R}, <)$ in a crucial way. If $X$ is a suborder of $\mathbb{R}$ containing a (complete) open interval $I$, then any partition $X = A \cup B$ of $X$ into two everywhere isomorphic sets would yield a partition of $I$ into two everywhere isomorphic sets. Since $I$ is isomorphic to $\mathbb{R}$, this would in turn yield a partition of $\mathbb{R}$ into two everywhere isomorphic sets, which is impossible by Theorem 1. The largest suborders of $\mathbb{R}$ that do not contain an open interval are those obtained by deleting some countable dense subset of $\mathbb{R}$. These are the largest suborders that might possibly be decomposed into two everywhere isomorphic sets. Since any such suborder is isomorphic to the irrationals $\mathbb{I} = \mathbb{R} \setminus \mathbb{Q}$, one may ask concretely if there is a partition $\mathbb{I} = A \cup B$ with $A$ and $B$ everywhere isomorphic. We show that, in contrast to $\mathbb{R}$, there is such a partition of $\mathbb{I}$. 

\theoremstyle{definition}
\newtheorem{thm2}[thm1]{Theorem}
\begin{thm2}\label{thm2} \,\ 
There exists a partition of the irrationals $\mathbb{I} = A \cup B$ such that $A$ and $B$ are everywhere isomorphic. 
\end{thm2}

The proof of Theorem \ref{thm1} also depends on the fact that we consider a partition of $\mathbb{R}$ into \emph{two} pieces. We do not know if it is possible to partition $\mathbb{R}$ into $n$ everywhere isomorphic sets for some $n > 2$. 

\theoremstyle{definition}
\newtheorem*{ques}{Question}
\begin{ques}\label{ques} \,\ 
Does there exist a partition $\mathbb{R} = A \cup B \cup C$ such that for every open interval $I$ we have $A \cap I \cong B \cap I \cong C \cap I$?
\end{ques}

Our notation and terminology are mostly standard. For linear orders $X$ and $Y$, we write $X \cong Y$ if $X$ is order-isomorphic to $Y$. An \emph{open interval} is a convex subset of a linear order with neither a top nor bottom point. If $X$ is a complete linear order (for example, if $X =\mathbb{R}$), every open interval of $X$ is of the form $(a, b), (a, \infty)$ $(-\infty, b)$, or $(-\infty, \infty)$ for some $a, b \in X$. But if $X$ has gaps, then it will have open intervals that cannot be written in terms of endpoints in $X$. A linear order $X$ is \emph{dense} if for any two distinct points in $X$ there is a third that lies strictly between them. A suborder $Y$ of a linear order $X$ is \emph{dense in $X$} if for any two distinct points of $X$, either both of them belong to $Y$ or there is a point between them that belongs to $Y$. Since $\mathbb{R}$ is dense as a linear order, any suborder $Y$ of $\mathbb{R}$ which is dense in $\mathbb{R}$ is also dense as a linear order.

\emph{Acknowledgments.} My warm thanks to Jackie Ferry, without whom this paper might never have quite been written; and to J\"urgen Kritschgau and Alexander Kechris for a number of useful comments and suggestions.

\section{Proof of Theorem 1}

Suppose toward a contradiction that there is a partition $\mathbb{R} = A \cup B$ such that for every open interval $I$, there exists an order-isomorphism $f_I: A \cap I \rightarrow B \cap I$. Then in particular $A$ and $B$ are each dense in $\mathbb{R}$, since they intersect every open interval. For a given open interval $I$, since $f_I$ is an order-preserving bijection \emph{from} a dense subset of $I$ \emph{onto} a dense subset of $I$, it can be extended uniquely to an order-automorphism of $I$. We identify $f_I$ with this extension. Then since $A \cap I$ and $B \cap I$ partition $I$ and $f_I[A \cap I] = B \cap I$, it must also be that $f_I[B \cap I] = A \cap I$. That is, $f_I$ interchanges the disjoint isomorphic suborders $A \cap I$ and $B \cap I$ while preserving the order of $I$. 

Since $f_I$ has no fixed points, given $x \in I$ we have that either $x < f_I(x)$ or $x > f_I(x)$. In the first case, since $f_I$ is order-preserving it must be that $f_I^n(x) < f_I^{n+1}(x)$ for all $n \in \mathbb{Z}$. That is we have
\[
\ldots < f_I^{-2}(x) < f_I^{-1}(x) < x < f_I(x) < f_I^2(x) < \ldots.
\]
Moreover, the increasing sequence of non-negative iterates $f_I^n(x), n \geq 0$, must be unbounded to the right in $I$. For if this sequence were bounded, so that the point $y = \lim_{n \rightarrow \infty} f_I^n(x)$ belonged to $I$, then by the continuity of $f_I$ we would have $f_I(y) = y$, a possibility we have ruled out. Likewise, the sequence of negative iterates $f_I^n(x), n < 0$, must be unbounded to the left in $I$. It follows, since $f_I$ is order-preserving, that it is increasing on all of $I$, that is $y < f_I(y)$ for all $y \in I$. 

On the other hand it may be that $x > f_I(x)$. In this case we have symmetrically
\[
\ldots < f_I^2(x) < f_I(x) < x < f_I^{-1}(x) < f_I^{-2}(x) < \ldots.
\]
The sequences of positive and negative iterates of $x$ are unbounded in $I$ to the left and right respectively, so that $f_I$ is decreasing on all of $I$. But then the inverse $f^{-1}_I$ is increasing on $I$, and this map is also an order-automorphism of $I$ that interchanges the sets $A \cap I$ and $B \cap I$. Thus by replacing $f_I$ with its inverse whenever necessary, we may assume that all of the $f_I$ are increasing. 

Observe that for any open interval $I$ and $x \in I$, the points in the iterate sequence
\[
\ldots < f_I^{-2}(x) < f_I^{-1}(x) < x < f_I(x) < f_I^2(x) < \ldots
\]
alternate in their belonging to $A$ or $B$: $f_I^{n}(x) \in A$ if and only if $f_I^{n+1}(x) \in B$. Thus if $J$ is an interval intersecting $I$ with $x \in I \cap J$, it can never be that an iterate of $x$ under $f_I$ is equal to an iterate of $x$ under $f_J$ of differing parity. Indeed, if for some $k, l \in \mathbb{Z}$ with $k \neq l \pmod 2$ there was a point $z$ such that $z = f_I^k(x) = f_J^l(x)$, then $z$ would belong to both $A$ and $B$, which is impossible. 

We will find such a $z$, and so obtain a contradiction. Let $K$ be a fixed open interval. We write $f$ for the automorphism $f_K$. Fix $x \in K$ and then fix a point $y$ in the interval $(f^{-1}(x), x)$. Let $I$ denote the open interval $(y, f^2(x))$. We write the automorphism $f_I$ as $g$. 

Since the increasing sequence of iterates $x < g(x) < g^2(x) < \ldots$ is unbounded to the right in $I$, there is a unique $n \geq 0$ such that $g^n(x) < f(x) \leq g^{n+1}(x)$. Let $J$ denote the closed interval $[g^n(x), f(x)]$. (It is convenient to consider closed intervals here, since we are going to apply the intermediate value theorem in a moment.) Then $J$ is an initial segment of the interval $[g^n(x), g^{n+1}(x)]$ and a final segment of the interval $[x, f(x)]$. Thus $f[J] = [f(g^n(x)), f^2(x)]$ is a final segment of $[f(x), f^2(x)]$. Since $f^2(x)$ is the right endpoint of $I$ we have that $f[J]$ is also a final segment of $I \cup \{f^2(x)\}$. For all sufficiently large $k$ we must have that $g^k(x) \in f[J]$, since these iterates are unbounded to the right in $I$. It follows that the intervals $g^k[J] = [g^{n+k}(x), g^k(f(x))]$ are subintervals of $f[J]$ for all sufficiently large $k$. 

Fix an even integer $N$ large enough so that $g^N[J]$ is a bounded subinterval of $f[J]$, that is, so that the left endpoint of $g^N[J]$ is strictly greater than the left endpoint of $f[J]$. We always have that the right endpoint of $g^N[J]$ is strictly less than the right endpoint of $f[J]$, since this endpoint is the same as the right endpoint of $I$. These maps are continuous on $J$ and hence so is their difference $f - g^N$. Since $f - g^N$ is negative at the left endpoint of $J$ and positive at the right endpoint, by the intermediate value theorem there must be an interior point $c \in J$ such that $(f - g^N)(c) = 0$, which gives $f(c) = g^N(c)$. Regardless of whether $c$ belongs to $A$ or $B$ we have that $f(c)$ belongs to the other set, whereas, by the evenness of $N$, $g^N(c)$ belongs to the same set as $c$. But $f(c) = g^N(c)$, contradicting $A \cap B = \emptyset$. \qed
\,\ \\

We note that the conclusion of Theorem \ref{thm1} holds more generally for any dense complete linear order $(L, <)$, so that if $A \cup B$ is a partition of such an order there must be an open interval on which the restrictions of $A$ and $B$ are not isomorphic. If $L$ does not embed $\omega_1$ or its reverse, this can be proved by essentially the same argument as for $\mathbb{R}$. If $L$ does embed either $\omega_1$ or its reverse, then the conclusion holds for a different reason. In such an $L$ there is always an open interval $I$ with an uncountable unbounded monotone sequence. Any automorphism of such an interval must have a fixed point, which prevents $I$ from being split into two isomorphic suborders both of which are dense in $I$. 

The completeness of $\mathbb{R}$ is used in several places in the proof, most significantly to guarantee that every iterate sequence $f_I^n(x)$ is unbounded in both directions in $I$, and again, in the guise of the intermediate value theorem, to produce a point in $A \cap B$. We cannot take corresponding steps if we partition a linear order that has gaps in every interval, and indeed, we show in the next section that the set of irrationals can be decomposed into two everywhere isomorphic sets. 

Nor can the proof be adapted to rule out the existence of a three-set partition $\mathbb{R} = A \cup B \cup C$ in which the restrictions $A \cap I, B \cap I$, and $C \cap I$ are pairwise isomorphic on every open interval $I$. What can be said is that there is no such partition in which the isomorphisms $A \cap I \cong B \cap I \cong C \cap I$ are simultaneously witnessed by a single automorphism $f_I: I \rightarrow I$, for every $I$. For example, if we always had $f_I[A \cap I] = B \cap I$ and $f_I[B \cap I] = C \cap I$ we could get an analogous contradiction. But this leaves open the possibility of such a partition for which, say, we have two isomorphisms $f_I: A \cap I \rightarrow B \cap I$ and $g_I: B \cap I \rightarrow C \cap I$ whose extensions to $I$ are distinct automorphisms of $I$, for every $I$. We conjecture there is such a partition. 

\section{Proof of Theorem 2}

We say that a suborder $X \subseteq \mathbb{R}$ is \emph{homogeneous} if for every open interval $I$ we have $X \cong X \cap I$.

For example, any countable dense subset $X$ of $\mathbb{R}$ is homogeneous, since by Cantor's theorem it holds that $X \cong X \cap I \cong \mathbb{Q}$ for every open interval $I$. In particular, $\mathbb{Q}$ is homogeneous. It turns out that for any infinite cardinal $\kappa \leq 2^{\aleph_0}$ there are homogeneous suborders $X \subseteq \mathbb{R}$ of size $\kappa$. In fact there are many such orders, and in the course of proving Theorem \ref{thm2} we will describe one method for constructing them. 

Our goal is to find a partition of the irrationals $\mathbb{I} = A \cup B$ into everywhere isomorphic sets $A$ and $B$. We will actually give such a decomposition in which both $A$ and $B$ are homogeneous. That is, for every open interval $I$ we will have $A \cong A \cap I \cong B \cap I \cong B$. (We emphasize that we will arrange this for \emph{every} open interval $I = (a, b)$ with $a, b \in \mathbb{R}$, not only for those with $a, b \in \mathbb{I}$.)

The reason for seeking a partition into homogeneous sets is that it will help us to reduce the number of isomorphisms we must find to prove isomorphism everywhere. For if $A$ is homogeneous, then there is a family of order-isomorphisms $f_I: A \rightarrow A \cap I$ witnessing this homogeneity. Identify these maps with their unique extensions to $\mathbb{R}$, so that $f_I: \mathbb{R} \rightarrow I$ is an order-isomorphism for every open interval $I$ with $f_I[A] = A \cap I$. If we also have $f_I[B] = B \cap I$ for every $I$, as we will arrange, then to show that $A$ and $B$ are everywhere isomorphic it suffices to find a single global isomorphism $g: A \rightarrow B$.  Then we obtain the isomorphisms $g_I: A \cap I \rightarrow B \cap I$, where $g_I = f_I \circ g \circ f_I^{-1}$, giving $A \cap I \cong B \cap I$ for every $I$.  

We will not work with the irrationals $\mathbb{I}$ directly but rather with the isomorphic linear order $\mathbb{Z}^{\omega}$. The reason for doing so is that the ordered group structure of $\mathbb{Z}^{\omega}$ will be used to define the global isomorphism $g$. Here $\omega = \{0, 1, 2, \ldots\}$ is the set of natural numbers, and $\mathbb{Z}^{\omega}$ is the set of sequences $u = (u_0, u_1, u_2, \ldots)$ with $u_i \in \mathbb{Z}$ for every $i \in \omega$. The ordering on $\mathbb{Z}^{\omega}$ is the lexicographical ordering: for distinct $u, v \in \mathbb{Z}^{\omega}$ we have $u < v$ if and only if $u_n < v_n$, where $n$ is the least integer in $\omega$ such that $u_n \neq v_n$. It is well-known that $\mathbb{Z}^{\omega}$ is order-isomorphic to $\mathbb{I}$, but since it will be helpful in what follows to see a proof, we sketch one.
\,\ \\

\underline{Claim}: $\mathbb{I} \cong \mathbb{Z}^{\omega}$.

\begin{proof}
For a dense linear order $X$ without endpoints, let $\overline{X}$ denote its Dedekind completion, the order obtained by filling each gap in $X$ with a single point. We do not add a top or bottom endpoint when taking the completion. 

Dedekind completions are complete, that is, every bounded monotone sequence in $\overline{X}$ converges to a point in $\overline{X}$. Up to isomorphism, $\overline{X}$ is the unique complete linear order (without endpoints) that contains $X$ as a dense suborder. Furthermore, for any $Y \subseteq \overline{X}$ which is dense in $\overline{X}$ we have $\overline{Y}=\overline{X}$. Thus if we can find an order $R$ (without endpoints) that is complete and contains $\mathbb{Z}^{\omega}$ as a dense suborder, it must be that $R \cong \overline{\mathbb{Z}^{\omega}}$. We define such an $R$ below. If moreover we can show that the difference $Q = R \setminus \mathbb{Z}^{\omega}$ is dense in $R$ and isomorphic to $\mathbb{Q}$, then it must be that $R = \overline{Q}$ is isomorphic to $\mathbb{R} = \overline{\mathbb{Q}}$, and $R \setminus Q = \mathbb{Z}^{\omega}$ is isomorphic to $\mathbb{R} \setminus \mathbb{Q} = \mathbb{I}$.  

Let $\mathbb{Z}^{<\omega}$ denote the set of nonempty finite sequences $r = (u_0, u_1, \ldots, u_n)$ with entries $u_i \in \mathbb{Z}$. Let $R = \mathbb{Z}^{\omega} \cup \mathbb{Z}^{<\omega}$. We order $R$ by the following rule: for distinct $u, v \in R$ we have $u < v$ if either there is an index $i$ for which $u_i \neq v_i$ and we have $u_i < v_i$ for the least such $i$, or $v$ is a finite sequence and $u$ extends $v$ as either a strictly longer finite sequence or infinite sequence. It is not hard to verify that this rule linearly orders $R$ and that both $\mathbb{Z}^{\omega}$ and $\mathbb{Z}^{<\omega}$ are dense in $R$. 

We claim furthermore that $R$ is complete. Since $\mathbb{Z}^{\omega}$ is dense in $R$, it is sufficient to check that every bounded monotone sequence in $\mathbb{Z}^{\omega}$ converges to a point in $R$. If $u^0 < u^1 < u^2 < \ldots$ is an increasing sequence in $\mathbb{Z}^{\omega}$, where $u^i = (u^i_0, u^i_1, \ldots)$, then either the first coordinates $u^i_0$ are unbounded in $\mathbb{Z}$, in which case the sequence is unbounded in $\mathbb{Z}^{\omega}$, or the first coordinates eventually stabilize, that is there is a $v_0 \in \mathbb{Z}$ such that for all sufficiently large $i$ we have $u^i_0 = v_0$. In this second case, there are likewise two possibilities. Either the second coordinates $u^i_1$ are unbounded in $\mathbb{Z}$, in which case the sequence converges to the point $(v_0) \in \mathbb{Z}^{<\omega}$, or the second coordinates also eventually stabilize, that is there is $v_1$ such that for all sufficiently large $i$ we have $u_1^i = v_1$. In this second case, we similarly have that either the sequence converges to $(v_0, v_1) \in \mathbb{Z}^{<\omega}$ or the third coordinates $u^i_2$ eventually stabilize at $v_2$. Continuing in this way, we either at some finite stage find a point $(v_0, v_1, \ldots, v_n) \in \mathbb{Z}^{<\omega}$ to which the sequence converges, or eventually all of its coordinates stabilize, in which case the sequence converges to a point $(v_0, v_1, \ldots) \in \mathbb{Z}^{\omega}$. Likewise it can be checked that bounded decreasing sequences must converge (though in this case if the sequence converges to a point in $\mathbb{Z}^{<\omega}$, it will be of the form $(v_0, v_1, \ldots, v_n-1)$). Hence $R$ is complete, as claimed. 

Since $\mathbb{Z}^{<\omega}$ is countable, has neither a top nor bottom point, and (being dense in $R$) is dense as a linear order, we have $\mathbb{Z}^{<\omega} \cong \mathbb{Q}$. It follows that $\overline{\mathbb{Z}^{<\omega}} \cong \mathbb{R}$ and $\overline{\mathbb{Z}^{<\omega}} \setminus \mathbb{Z}^{<\omega} \cong \mathbb{I}$. But $\overline{\mathbb{Z}^{<\omega}} = R$ by above and $R \setminus \mathbb{Z}^{<\omega} = \mathbb{Z}^{\omega}$, giving $\mathbb{Z}^{\omega} \cong \mathbb{I}$. 
\end{proof}

We will construct the everywhere isomorphic suborders $A$ and $B$ in the statement of Theorem \ref{thm2} as disjoint suborders of $\mathbb{Z}^{\omega}$. At the end, we will not have arranged that $A \cup B = \mathbb{Z}^{\omega}$, but rather $A \cup B = \mathbb{Z}^{\omega} \setminus C$, where $C$ is some countable subset of $\mathbb{Z}^{\omega}$. But the order type of the irrationals is not changed by deleting a countable set of points, since after such a deletion its complement in the reals, which has gained countably many points, nonetheless remains countable and dense and therefore isomorphic to $\mathbb{Q}$. Hence we will have succeeded in getting the desired decomposition, since $A \cup B = \mathbb{Z}^{\omega} \setminus C \cong \mathbb{I} \setminus C \cong \mathbb{I}$.

We introduce notation for dealing with sequences. We use $u, v, \ldots$ for infinite sequences (elements of $\mathbb{Z}^{\omega}$), $r, s, \ldots$ for finite sequences (elements of $\mathbb{Z}^{<\omega}$), and $n, m, \ldots$ for elements of $\mathbb{Z}$. We sometimes view elements of $\mathbb{Z}$ as sequences of length 1, and think of $\mathbb{Z}$ as a subset of $\mathbb{Z}^{<\omega}$. If $r$ is a finite sequence and $u$ is either a finite or infinite sequence, we write $ru$ for the sequence of $r$ concatenated with $u$. For a sequence $v$ written as a concatenation $v = ru$, we say that $r$ is an \emph{initial sequence} of $v$, and $u$ is a \emph{tail-sequence}.

For sequences $u, v \in \mathbb{Z}^{\omega}$ we say $u$ and $v$ are \emph{tail-equivalent} and write $u \sim v$ if there exist finite sequences $r, s \in \mathbb{Z}^{<\omega}$ and an infinite sequence $u' \in \mathbb{Z}^{\omega}$ such that $u = ru'$ and $v = su'$. We emphasize that $r$ and $s$ need not be of the same length in this definition. Tail-equivalence is an equivalence relation on $\mathbb{Z}^{\omega}$. We denote the equivalence class of a given $u \in \mathbb{Z}^{\omega}$ by $[u]$. It can be checked that $[u]$ is exactly the set of sequences of the form $ru'$, where $r \in \mathbb{Z}^{<\omega}$ is an arbitrary finite sequence and $u'$ is a tail-sequence of $u$. Since there are only countably many tails $u'$ of $u$ and countably many $r \in \mathbb{Z}^{<\omega}$, it follows that each tail-equivalence class is countable. Since these classes partition $\mathbb{R}$ into countable sets, the number of classes is $2^{\aleph_0}$.

If $u$ and $v$ are tail-equivalent, we call a decomposition $u = ru'$ and $v = su'$ a \emph{meeting representation} of $u$ and $v$. Meeting representations are not unique, since we can always further unzip along the tail-sequence $u'$ to get new representations. That is, if we write $u' = tu''$ for some initial sequence $t$ of $u'$, then $u = rtu''$ and $v = stu''$ is also a meeting representation of $u$ and $v$, now with respect to the tail-sequence $u''$ instead of $u'$. 

For $r \in \mathbb{Z}^{<\omega}$, let $I_r$ denote the set of sequences in $\mathbb{Z}^{\omega}$ beginning with $r$, so that $u \in I_r$ if and only if there is $u' \in \mathbb{Z}^{\omega}$ such that $u=ru'$. As a subset of $\mathbb{Z}^{\omega}$, each $I_r$ is an open interval, that we call a \emph{standard interval}. The endpoints of $I_r$ do not lie in $\mathbb{Z}^{\omega}$ but rather in its completion $R = \mathbb{Z}^{\omega} \cup \mathbb{Z}^{<\omega}$: if $r = (u_0, u_1, \ldots, u_n)$ and we let $r' = (u_0, u_1, \ldots, u_n-1)$, then $I_r$ is exactly the open interval $(r', r)$. More strictly speaking, since we are viewing $I_r$ as a subset of $\mathbb{Z}^{\omega}$ only, it is the set of points in $\mathbb{Z}^{\omega}$ lying between $r'$ and $r$. Observe that for two finite sequences $r$ and $s$, if $r$ is lexicographically less than $s$, then $I_r$ lies entirely to the left of $I_s$, whereas if $r$ extends $s$ then $I_r$ is a subinterval of $I_s$. 

Let $f_r$ denote the \emph{projection} onto $I_r$, defined by $f_r(u) = ru$ for all $u \in \mathbb{Z}^{\omega}$. It is quick to check that $f_r$ is an order-isomorphism of $\mathbb{Z}^{\omega}$ with $I_r$. Since we have $f_r \circ f_s = f_{rs}$ for all $r, s \in \mathbb{Z}^{<\omega}$, we have that all of the $f_r$ can be written as iterated compositions of the projections $f_n$ for $n \in \mathbb{Z}$.

The projections $f_r$ are closely related to the tail-equivalence relation. We claim that for a fixed $r \in \mathbb{Z}^{<\omega}$ and $u \in \mathbb{Z}^{\omega}$ we have that $f_r[[u]] = [u] \cap I_r$. That is, each $f_r$ witnesses that the tail-equivalence class $[u]$ is order-isomorphic to its restriction to $I_r$. If $v \in [u]$, so that $v = su'$ for some finite sequence $s$ and some tail-sequence $u'$ of $u$, then $f_r[v] = rv = rsu'$ is also tail-equivalent to $u$, giving the forward containment. On the other hand, if $w \in [u] \cap I_r$, observe that $w$ can be written $w = rsu'$ for some finite sequence $s$ and tail-sequence $u'$ of $u$. (To write $w$ as such we may need to unzip along the tail-sequence of a given meeting representation of $u$ and $w$, but this is no problem.) But then we have $f_r[su'] = w$, and since $su' \in [u]$ this gives the reverse containment. It follows that if $X$ is a suborder of $\mathbb{Z}^{\omega}$ that is closed under tail-equivalence (that is, $X$ is a union of tail-equivalence classes), then $f_r[X] = X \cap I_r$ for every $r \in \mathbb{Z}^{<\omega}$.

We note that it holds conversely that if $f_r[X] = X \cap I_r$ for every $r$, then $X$ is closed under tail-equivalence. This is not hard to check, and we leave it to the reader since we will not need it for our construction. Thus the tail-equivalence classes are the smallest suborders of $\mathbb{Z}^{\omega}$ that are invariant under all of the $f_r$.  

Our aim is to construct suborders $A$ and $B$ of $\mathbb{Z}^{\omega}$ that, in addition to being everywhere isomorphic, are homogeneous. To get homogeneity, we will arrange that $A$ and $B$ are unions of tail-equivalence classes. By what we have just shown, this guarantees that $A \cong A \cap I$ and $B \cong B \cap I$ for every interval of the form $I = I_r$. It turns out that this is enough to get $A \cong A \cap I$ and $B \cong B \cap I$ for \emph{every} open interval $I$. 
\,\ \\

\underline{Claim}: Suppose $X \subseteq \mathbb{Z}^{\omega}$ is closed under tail-equivalence, so that $X \cong X \cap I_r$ for every $r \in \mathbb{Z}^{<\omega}$. Then $X$ is homogeneous.
\,\ \\

In proving the claim we will be dealing with segments of orders and sums of orders. An \emph{initial segment} $I$ of a linear order $X$ is an interval in $X$ that is unbounded to the left, that is, if $y \in I$ and $x < y$ then $x \in I$. A \emph{final segment} $J$ is the complement of an initial segment, or equivalently, an interval that is unbounded to the right. 

We will need the following fact, due to Lindenbaum (\cite{Lindenbaum}, or see \cite[pg. 248]{Sierpinski}): if $L$ and $M$ are linear orders such that $L$ is isomorphic to an initial segment of $M$ and $M$ is isomorphic to a final segment of $L$, then $L$ is isomorphic to $M$. We will prove a refined version of Lindenbaum's theorem below.

Given linear orders $L_0$ and $L_1$, we write $L_0 + L_1$ for the order, unique up to isomorphism, that can be partitioned into an initial segment isomorphic to $L_0$ and corresponding final segment isomorphic to $L_1$. We also consider longer sums $L_0 + L_1 + \ldots + L_n$, or even infinite sums, such as \emph{$\omega$-sums} 
\[
L_0 + L_1 + L_2 + \ldots,
\] 
or \emph{$\omega^*$-sums} 
\[
\ldots + L_2 + L_1 + L_0,
\] 
or \emph{$\mathbb{Z}$-sums} 
\[
\ldots + L_{-1} + L_0 + L_1 + L_2 + \ldots.
\]
If all of the orders $L_i$ are isomorphic to a single order $L$, then we write $nL$ for the $n$-sum, $\omega L$ for the $\omega$-sum, $\omega^*L$ for the $\omega^*$-sum, and $\mathbb{Z}L$ for the $\mathbb{Z}$-sum. 

Formally here, we are viewing a product of two orders $XY$ as the cartesian product $X \times Y$ of their underlying sets, ordered lexicographically by the rule $(x, y) < (x', y')$ if either $x < x'$, or $x = x'$ and $y < y'$. So, for example, $\omega L = \{(n, l): n \in \omega, l \in L\}$, ordered lexicographically.
\,\ \\ \,\ \\
\emph{Proof of claim.}
For $r \in \mathbb{Z}^{<\omega}$ we write $X_r$ for $X \cap I_r$. The orders $I_n$, $n \in \mathbb{Z}$, partition $\mathbb{Z}^{\omega}$ into $\mathbb{Z}$-many copies of itself. Since $X \cong X_n$ for every $n$, it follows that the $X_n$ partition $X$ into $\mathbb{Z}$-many copies of itself. That is, we have $X \cong \mathbb{Z}X$.

We claim that also $X \cong \omega X$. On one hand, $X$ is isomorphic to an initial segment of $\omega X$, namely the initial copy of itself in the $\omega$-sum. On the other, $\omega X$ is naturally isomorphic to the right half of $\mathbb{Z}X$, which is a final segment of $\mathbb{Z}X$. But $\mathbb{Z}X$ is isomorphic to $X$, so $\omega X$ is isomorphic to a final segment of $X$. By Lindenbaum's theorem we have that $X \cong \omega X$ as claimed.

We can similarly prove that $X \cong \omega^*X$. Then, by splitting $\mathbb{Z}X$ into its left and right halves, we have $\mathbb{Z}X \cong \omega^*X + \omega X$, which gives $\mathbb{Z}X \cong 2X$, and hence $X \cong 2X$. Repeated application of this last identity gives $X \cong nX$ for $n \geq 3$ as well. 

The identities $X \cong \mathbb{Z}X \cong \omega X \cong \omega^*X \cong nX$ make it possible to prove that $X \cong X \cap I$ for any open interval $I$, by showing that $X \cap I$ can be written as a finite sum of terms each of which is isomorphic to one of these orders. To do this, we analyze the form of open intervals in $\mathbb{Z}^{\omega}$ according to their endpoints, remembering that these endpoints may lie in $\mathbb{Z}^{<\omega}$. We analyze unbounded intervals first, since these are determined by a single endpoint. 

Suppose that $I$ is an open interval that is unbounded to the right in $\mathbb{Z}^{\omega}$, so that $I$ is a final segment of $\mathbb{Z}^{\omega}$. We assume that $I$ is not also an initial segment, that is, $I \neq \mathbb{Z}^{\omega}$. There are two possibilities. Either $I = (u, \infty)$ for some $u \in \mathbb{Z}^{\omega}$ or $I = (r, \infty)$ for some $r \in \mathbb{Z}^{<\omega}$. Suppose we are in this second case, and $r = (u_0, u_1, \ldots, u_n)$. Then $I$ consists of the infinite sequences that are lexicographically greater than $r$, those $v \in \mathbb{Z}^{\omega}$ for which there is $k \leq n$ such that $v_k \neq u_k$, and for the least such $k$ we have $v_k > u_k$. Such sequences are categorized by which coordinate $k \leq n$ they first differ from $r$. If $v_0 > u_0$ then $v$ belongs to one of the intervals $I_m$, where $m \in \mathbb{Z}$ and $m > u_0$. The collection of these intervals is the $\omega$-sum
\[
\sum_{i = 1}^{\infty} I_{u_0 + i} = I_{u_0+1} + I_{u_0+2} + \ldots.
\]
This sum is a final segment of $I$ (and also of $\mathbb{Z}^{\omega}$). 

If instead $v_0 = u_0$, so that $v$ belongs to $I_{u_0}$, but we have $v_1 > u_1$, then $v$ belongs to one of the intervals $I_{(u_0, m)}$ where $m > u_1$. These form the $\omega$-sum

\[
\sum_{i = 1}^{\infty} I_{(u_0, u_1 + i)} = I_{(u_0, u_1 + 1)} + I_{(u_0, u_1 + 1)} + \ldots
\]
which is a final segment of $I_{u_0}$. Continuing in this way, we see that $v$ belongs to one of the intervals $I_{(u_0, \ldots, u_{k-1}, m)}$ where $k \leq n$ and $m > u_k$. For a fixed $k \leq n$ these intervals form an $\omega$-sum, and each of these sums lies immediately to the left of the previous one, so that we can decompose $I$ as follows:
\begin{align*}
I & = \sum_{i = 1}^{\infty} I_{(u_0, u_1, \ldots, u_n + i)} + \ldots + \sum_{i = 1}^{\infty} I_{(u_0, u_1 + i)} + \sum_{i = 1}^{\infty} I_{u_0 + i} \\
& = \sum_{j = n}^0 \sum_{i = 1}^{\infty} I_{(u_0, u_1, \ldots, u_j + i)}.
\end{align*}
This represents $I$ as an $(n+1)$-sum of $\omega$-sums of standard intervals. It follows that 
\begin{align*}
X \cap I & = \sum_{j = n}^0 \sum_{i = 1}^{\infty} X \cap I_{(u_0, u_1, \ldots, u_j + i)}.
\end{align*}
But $X$ is isomorphic to its restriction to every standard interval, so we have
\[
X \cap I  \cong \sum_{j = n}^0 \sum_{i = 1}^{\infty} X = \sum_{j = n}^0 \omega X \cong \sum_{j = n}^0 X = (n+1)X \cong X,
\]
as desired.

Suppose now that $I = (u, \infty)$ for some point $u = (u_0, u_1, \ldots)$ in $\mathbb{Z}^{\omega}$. Going coordinate by coordinate as above, we may decompose $I$ similarly, not as a finite sum of $\omega$-sums of standard intervals, but an $\omega^*$-sum of $\omega$-sums of standard intervals:
\[
I = \ldots + \sum_{i = 1}^{\infty} I_{(u_0, u_1 + i)} + \sum_{i = 1}^{\infty} I_{u_0 + i}.
\]
Abusing notation, we write
\[
I = \sum_{j = \infty}^0 \sum_{i = 1}^{\infty} I_{(u_0, u_1, \ldots, u_j + i)}.
\]
Then we have
\[
X \cap I  \cong \sum_{j = \infty}^0 \sum_{i = 1}^{\infty} X = \sum_{j = \infty}^0 \omega X \cong \sum_{j = \infty}^0 X = \omega^*X \cong X.
\]

Now suppose that $I$ is unbounded to the left, so that $I$ is an initial segment of $\mathbb{Z}^{\omega}$. Then either $I = (-\infty, u)$ for some $u \in \mathbb{Z}^{\omega}$ or $I = (-\infty, r)$ for some $r \in \mathbb{Z}^{<\omega}$. Both cases mirror the corresponding cases above. If $I = (-\infty, r)$, then $I$ can be decomposed as a finite sum of $\omega^*$-sums of standard intervals, so that $X \cap I$ can be decomposed as a finite sum of $\omega^*$-sums of copies of $X$, which is isomorphic to $X$. If $I = (-\infty, u)$ then $I$ can be decomposed as an $\omega$-sum of $\omega^*$-sums of standard intervals, so that $X \cap I$ can be decomposed as an $\omega$-sum of $\omega^*$-sums of copies of $X$, which is also isomorphic to $X$. 

Finally it may be that $I$ is bounded on both sides, so that $I = (x, y)$ for some $x, y \in \mathbb{Z}^{\omega} \cup \mathbb{Z}^{<\omega}$. We just showed that for any $z \in \mathbb{Z}^{\omega} \cup \mathbb{Z}^{<\omega}$ we have $X \cap (-\infty, z) \cong X$. Another way of expressing this is that any initial segment of $X$ without a top point is isomorphic to $X$. But $X \cap I$ is an initial segment of $X \cap (x, \infty)$ without a top point, and this latter order we know to be isomorphic to $X$. Hence we have that $X \cap I$ is isomorphic to $X$. \qed
\,\ \\

Observe that it follows from the claim that for any infinite cardinal $\kappa \leq 2^{\aleph_0}$ there exists a homogeneous order $X \subseteq \mathbb{R}$ of size $\kappa$: take $X$ to be the union of any $\kappa$-many distinct tail-equivalence classes. Such an $X$ is constructed as a homogeneous suborder of $\mathbb{Z}^{\omega}$, but since $\mathbb{Z}^{\omega}$ is isomorphic to the irrationals, $X$ is isomorphic to a homogeneous suborder of $\mathbb{R}$. 

Before we can define the sets $A$ and $B$ from Theorem 2, we need to refine the claim above to handle the case when we are dealing with not one, but rather two disjoint suborders $A, B$ of $\mathbb{Z}^{\omega}$ each of which is closed under tail-equivalence. We wish to conclude in such a situation that not only do we have for each open interval $I$ that $A \cong A \cap I$ and $B \cong B \cap I$, but actually that there is a single isomorphism $f_I: \mathbb{Z}^{\omega} \rightarrow I$ such that $f_I[A] = A \cap I$ and $f_I[B] = B \cap I$. As we observed, this reduces proving that $A \cap I \cong B \cap I$ on every $I$ to finding a single isomorphism $g: A \rightarrow B$. That such isomorphisms $f_I$ exist is implicit in the proof of the claim above. We draw them out explicitly. 

We need the following refinement of Lindenbaum's theorem.
\,\ \\

\underline{Lemma}: Suppose $X$ and $Y$ are linear orders, and $X_0 \subseteq X$ and $Y_0 \subseteq Y$ are suborders of $X$ and $Y$ respectively. Suppose $f: X \rightarrow Y$ is an embedding of $X$ onto an initial segment of $Y$ such that $f[X_0] = Y_0 \cap f[X]$, and $g: Y \rightarrow X$ is an embedding of $Y$ onto a final segment of $X$ such that $g[Y_0] = X_0 \cap g[Y]$. Then there is an isomorphism $h: X \rightarrow Y$ such that $h[X_0] = Y_0$.

\begin{proof}
Let $h$ be the bijection built out of $f$ and $g$ as in the standard proof of the Schroeder-Bernstein theorem. Observe that the hypotheses that $f$ is onto an initial segment of $Y$ and $g$ is onto a final segment of $X$ guarantee that $h$ is order-preserving, that is, an isomorphism of $X$ with $Y$. Observe moreover that since the same hypotheses apply to the restrictions of $f$ to $X_0$ and $g$ to $Y_0$, we have that $h[X_0] = Y_0$. 
\end{proof}

Now, suppose that $X = \mathbb{Z}^{\omega}$. Applying the claim above gives that $\mathbb{Z}^{\omega} \cong \mathbb{Z}^{\omega} \cap I$ for any open interval $I$. (Since all of our maps are into $\mathbb{Z}^{\omega}$, we will usually write $\mathbb{Z}^{\omega} \cap I$ simply as $I$ below.) But the proof of the claim implicitly shows more. The natural isomorphisms that witness the identities $\mathbb{Z}^{\omega} \cong \mathbb{Z}\mathbb{Z}^{\omega} \cong \omega\mathbb{Z}^{\omega} \cong \omega^*\mathbb{Z}^{\omega} \cong n\mathbb{Z}^{\omega}$ are all combinations of the projection maps $f_r$. Since these maps preserve tail-equivalence, in the strong sense that $f_r[A] = A \cap I_r$ for any $A \subseteq \mathbb{Z}^{\omega}$ that is closed under tail-equivalence, it follows that the isomorphisms $f_I: \mathbb{Z}^{\omega} \rightarrow I$ we get out of these identities will have the property that $f_I[A] = A \cap I$ for any such $A$.

Let us show this explicitly. Fix a single tail-equivalence $A$ (so that $A = [u]$ for any $u \in A$). The identity $\mathbb{Z}\mathbb{Z}^{\omega} \cong \mathbb{Z}^{\omega}$ is naturally witnessed by the flattening map $fl: \mathbb{Z}\mathbb{Z}^{\omega} \rightarrow \mathbb{Z}^{\omega}$ defined by $fl(z, u) = zu$. It is easily checked that this is an order-isomorphism of $\mathbb{Z}\mathbb{Z}^{\omega}$ and $\mathbb{Z}^{\omega}$, and moreover that we have $fl[\mathbb{Z}A] = A$, where $\mathbb{Z}A = \{(z, u) \in \mathbb{Z}\mathbb{Z}^{\omega}: u \in A\}$. 

The initial copy of $\mathbb{Z}^{\omega}$ in the order $\omega \mathbb{Z}^{\omega}$ is the interval $I = \{(0, u): u \in \mathbb{Z}^{\omega}\}$ consisting of points whose first coordinate is $0$. (This $I$ is essentially the interval $I_0$.) That $\mathbb{Z}^{\omega}$ is isomorphic to this initial copy of itself is witnessed by the map $f: \mathbb{Z}^{\omega} \rightarrow \omega \mathbb{Z}^{\omega}$ defined by $f(u) = (0, u)$. (This $f$ is essentially the projection map $f_0$.) We clearly have that $f[A] = \omega A \cap I = \omega A \cap f[\mathbb{Z}^{\omega}]$. On the other hand, the natural map witnessing that $\omega \mathbb{Z}^{\omega}$ is isomorphic to a final segment of $\mathbb{Z}^{\omega}$ is just the flattening map $fl$, restricted to $\omega \mathbb{Z}^{\omega}$. It is quickly checked that $fl[\omega A] = A \cap fl[\omega \mathbb{Z}^{\omega}]$. The lemma then yields an isomorphism $f_{\omega}: \mathbb{Z}^{\omega} \rightarrow \omega \mathbb{Z}^{\omega}$ that sends $A$ onto $\omega A$. 

By a symmetric argument we get an isomorphism $f_{\omega^*}: \mathbb{Z}^{\omega} \rightarrow \omega^* \mathbb{Z}^{\omega}$ that maps $A$ onto $\omega^*A$. Now, viewing $\mathbb{Z}\mathbb{Z}^{\omega}$ as being composed of the initial segment $\omega^*\mathbb{Z}^{\omega}$ followed by final segment $\omega\mathbb{Z}^{\omega}$, we get an isomorphism from $\mathbb{Z}\mathbb{Z}^{\omega}$ to $2\mathbb{Z}^{\omega}$ by the rule $(k, u) \mapsto (0, f_{\omega^*}^{-1}((k,u))$ if $k < 0$ and $(k, u) \mapsto (1, f_{\omega}^{-1}((k,u))$ if $k \geq 0$. (Here, we are identifying 2 with $\{0, 1\}$ and, for convenience, $\omega^*$ with $\{\ldots, -3, -2, -1\}$.) This map sends $\mathbb{Z}A$ onto $2A$. By composing it with $fl^{-1}$, we get an isomorphism $f_2: \mathbb{Z}^{\omega} \rightarrow 2\mathbb{Z}^{\omega}$ that sends $A$ onto $2A$. To get a map onto $3\mathbb{Z}^{\omega}$, we may apply $f_2$ to the righthand copy of $\mathbb{Z}^{\omega}$ in $2\mathbb{Z}^{\omega}$. Explicitly, let $3'\mathbb{Z}^{\omega}$ denote the set of all tuples of the form $(0, u)$, $(1, (0, u))$, and $(1, (1, u))$ with $u \in \mathbb{Z}^{\omega}$. View $3'\mathbb{Z}^{\omega}$ as being ordered lexicographically in the natural way. Define an isomorphism $g_{3'}: 2\mathbb{Z}^{\omega} \rightarrow 3'\mathbb{Z}^{\omega}$ by the rule $g_{3'}((0,u)) = (0,u)$ and $g_{3'}((1, u)) = (1, f_2(u))$.  Compose $g_{3'}$ with the map defined by the rules $(0, u) \mapsto (0, u); (1, (0, u)) \mapsto (1, u); (1, (1, u)) \mapsto (2, u)$ to get an isomorphism $g_3$ from $2\mathbb{Z}^{\omega}$ onto $3\mathbb{Z}^{\omega}$. (Here, $3 = \{0, 1, 2\}$.) Finally, let $f_3 = g_3 \circ f_2$. Then $f_3: \mathbb{Z}^{\omega} \rightarrow 3 \mathbb{Z}^{\omega}$ is an isomorphism that takes $A$ onto $3A$. Similarly we get isomorphisms $f_n: \mathbb{Z}^{\omega} \rightarrow n \mathbb{Z}^{\omega}$ that take $A$ onto $nA$, for every $n \in \omega$. 

Thus we have shown that our maps witnessing $\mathbb{Z}^{\omega} \cong \mathbb{Z}\mathbb{Z}^{\omega} \cong \omega \mathbb{Z}^{\omega} \cong \omega^*\mathbb{Z}^{\omega} \cong n \mathbb{Z}^{\omega}$ also witness $A \cong \mathbb{Z} A \cong \omega A \cong \omega^* A \cong nA$. From these maps, we can build isomorphisms $f_I: \mathbb{Z}^{\omega} \rightarrow I$ that send $A$ onto $A \cap I$, for each open interval $I$. 

For instance, when $I = (r, \infty)$ for some $r \in \mathbb{Z}^{<\omega}$, we decompose $I$ as in the proof of the claim:
\begin{align*}
I & = \sum_{j = n}^0 \sum_{i = 1}^{\infty} I_{(u_0, u_1, \ldots, u_j + i)}.
\end{align*}
For a fixed $j$, $0 \leq j \leq n$, and $i$, $1 \leq i < \infty$, let $r(j, i)$ denote the finite sequence $(u_0, \ldots, u_j+i)$. Each summand $I_{(u_0, \ldots, u_j+i)} = I_{r(j, i)}$ is naturally isomorphic to $\mathbb{Z}^{\omega}$ via the map $f_{r(j, i)}^{-1}$. The sum 
\[
\sum_{i = 1}^{\infty} I_{(u_0, u_1, \ldots, u_j + i)}
\]
is isomorphic to $\omega \mathbb{Z}^{\omega}$ via the map $g$ defined by the following rule: for each $i$, $1 \leq i < \infty$, if $u \in I_{r(j, i)}$ then $g(u) = (i-1, f_{r(j, i)}^{-1}(u))$. (We shift to $i-1$ merely to be correct, since $\omega$ begins at $0$.) Since $f_{r(j, i)}^{-1}$ maps $A \cap I_{r(j, i)}$ onto $A$, we have 
\[
g\left[A \cap \sum_{i = 1}^{\infty} I_{(u_0, u_1, \ldots, u_j + i)}\right] = \omega A.
\]
Then, since our maps witnessing $\omega \mathbb{Z}^{\omega} \cong \mathbb{Z}^{\omega}$ and $(n+1) \mathbb{Z}^{\omega} \cong \mathbb{Z}^{\omega}$ send $\omega A$ onto $A$ and $(n+1)A$ onto $A$ respectively, we have that the isomorphisms witnessing 
\[
I = \sum_{j = n}^0 \sum_{i = 1}^{\infty} I_{(u_0, u_1, \ldots, u_j + i)} \cong \sum_{j = n}^0 \omega \mathbb{Z}^{\omega} \cong \sum_{j = n}^0 \mathbb{Z}^{\omega} = (n+1)\mathbb{Z}^{\omega} \cong \mathbb{Z}^{\omega}
\]
also witness
\[
A \cap I = \sum_{j = n}^0 \sum_{i = 1}^{\infty} A \cap I_{(u_0, u_1, \ldots, u_j + i)} \cong \sum_{j = n}^0 \omega A \cong \sum_{j = n}^0 A = (n+1)A \cong A.
\]
Thus there is an isomorphism $f_I: \mathbb{Z}^{\omega} \rightarrow I$ such that $f_I[A] = A \cap I$.

A similar argument goes through when the left endpoint of $I$ is some $u \in \mathbb{Z}^{\omega}$, and also in the two symmetric cases when $I$ is unbounded to the left. 

So suppose we are in the last case, when $I = (x, y)$ for points $x, y \in \mathbb{Z}^{\omega} \cup \mathbb{Z}^{<\omega}$. Let $J = (x, \infty)$. We know we have an isomorphism $f_J: \mathbb{Z}^{\omega} \rightarrow J$ such that $f[A] = A \cap J$. Let $I' = f_J^{-1}[I]$. Since $I$ is an initial segment of $J$, we have that $I'$ is an initial segment of $\mathbb{Z}^{\omega}$, and $f_J^{-1}[A \cap I] = A \cap I'$. We also have an isomorphism $f_I': \mathbb{Z}^{\omega} \rightarrow I'$ such that $f_I'[A] = A \cap I'$. But then the map $f_I = f_J \circ f_I'$ is an isomorphism of $\mathbb{Z}^{\omega}$ with $I$ such that $f_I[A] = A \cap I$. 

Thus for every open interval $I$ we have an isomorphism $f_I: \mathbb{Z}^{\omega} \rightarrow I$ such that $f_I[A] = A \cap I$. Since $A$ was an arbitrary tail-equivalence class and our construction of the $f_I$ did not depend on the particular class we fixed, the same statement holds when $A$ is any union of tail-equivalence classes. We have proved the following refinement of our previous claim.
\,\ \\

\underline{Claim}: For every open interval $I \subseteq \mathbb{Z}^{\omega}$, there is an isomorphism $f_I: \mathbb{Z}^{\omega} \rightarrow I$ such that for every suborder $A \subseteq \mathbb{Z}^{\omega}$ that is closed under tail-equivalence we have $f_I[A] = A \cap I$.
\,\ \\

We are nearly ready to define the sets $A$ and $B$ from Theorem 2. Both $A$ and $B$ will be unions of tail-equivalence classes. By our refined claim, this guarantees that for every open $I$ we have an isomorphism $f_I: \mathbb{Z}^{\omega} \rightarrow I$ that takes both $A$ onto $A \cap I$ and $B$ onto $B \cap I$. Thus to ensure that $A$ and $B$ are everywhere isomorphic, we need only show that there is a global isomorphism $g: A \rightarrow B$. Our $g$ will actually be an order-automorphism of $\mathbb{Z}^{\omega}$ such that $g[A] = B$ and $g[B] = A$. We define $g$ first, and then construct $A$ and $B$ to satisfy these identities while still being closed under tail-equivalence. 

To define $g$, we avail ourselves of the ordered group structure of $\mathbb{Z}^{\omega}$. For sequences $u = (u_0, u_1, \ldots)$ and $v = (v_0, v_1, \ldots)$ we write $u + v$ for the sequence $(u_0 + v_0, u_1 + v_1, \ldots)$. For any fixed $v \in \mathbb{Z}^{\omega}$ the map $u \mapsto u + v$ defines an order-automorphism of $\mathbb{Z}^{\omega}$. We choose a specific one. For $n \in \mathbb{Z}$, write $\overline{n}$ for the sequence $(n, n, n, \ldots)$. Define $g: \mathbb{Z}^{\omega} \rightarrow \mathbb{Z}^{\omega}$ by the rule $g(u) = u + \overline{1}$. 

That $g$ is a shift by $\overline{1}$ in particular is not important, except that $\overline{1}$ is not tail-equivalent to the identity $\overline{0}$. For the construction to work, we need that $g(u) \not\sim u$ for all but countably many $u$. 

Observe that $g$ preserves tail-equivalence in the weaker sense that $u \sim v$ if and only if $g(u) \sim g(v)$. Thus $g$ induces a permutation of the tail-equivalence classes. Specifically, for every $u \in \mathbb{Z}^{\omega}$ we have $g[[u]] = [g(u)] = [u + \overline{1}]$. 

For a fixed $u$, consider the iterated images of the tail-equivalence class of $u$ under $g$. These are the sets $g^k[[u]] = [u+\overline{k}]$ for $ k \in \mathbb{Z}$. We call the union $\bigcup_{k \in \mathbb{Z}} [u + \overline{k}]$ the \emph{orbit of} $u$ and denote it by $O(u)$. Note that every orbit is countable. 

For most $u$ we have $[u + \overline{k}] \cap [u + \overline{l}] = \emptyset$ whenever $k \neq l$. But not for all $u$. Let $C$ denote the union of all orbits $O(u)$ for which the classes $[u + \overline{k}], k \in \mathbb{Z}$ are not pairwise disjoint. It is not hard to see that $C$ is exactly the set of $u$ for which there exists $k \neq 0$ such that $u \sim u + \overline{k}$. 

We claim that $C$ is countable. Fix $u \in C$, and find $k \neq 0$ such that $u \sim u + \overline{k}$. Then $u = ru'$ and $u + \overline{k} = su'$ for some $r, s \in \mathbb{Z}^{<\omega}$ and $u' \in \mathbb{Z}^{\omega}$. We write $|r|$ and $|s|$ for the lengths of the sequences $r$ and  $s$ respectively. It cannot be that $|r| = |s|$, since if this were so, the tail-sequence $u'$ would begin at the same coordinate in both $u$ and $u + \overline{k}$, giving the false identity $u' = u' + \overline{k}$. So suppose that $|r| > |s|$. The case when $|r| < |s|$ is symmetric. Write $r$ as $r't$ where $|r'| = |s|$ and $t$ is the remainder sequence. Then we have $u = r'tu'$ and $u + \overline{k} = su'$. Let us decompose $u'$ into finite blocks $q_i$ by writing $u' = q_0q_1q_2\ldots$, where each block $q_i$ has length $|t|$. Then we have $u = r'tq_0q_1\ldots$ and $u + \overline{k}=sq_0q_1\ldots$. For a finite block $q$ we write $q + \overline{k}$ for the finite block of the same length as $q$ in which $k$ has been added to every entry. We have decomposed $u$ and $u + \overline{k}$ so that each block in $u$ has the same length as the corresponding block in $u + \overline{k}$. Thus we obtain the equations $s = r' + \overline{k}$, $q_0 = t + \overline{k}$, $q_1 = q_0 + \overline{k}$, $q_2 = q_1 + \overline{k}$, etc. But then $u = r't(t+\overline{k})(t+\overline{2k})\ldots$. Such a $u$ is specified by the initial sequence $r$, the block $t$, and the difference $k$. There are only countably many $r, t \in \mathbb{Z}^{<\omega}$ and countably many $k \in \mathbb{Z}$. Hence $C$ is countable, as claimed. 

We construct $A$ and $B$. Consider $\mathbb{Z}^{\omega} \setminus C$. This set consists of the orbits $O(u)$ not belonging to $C$, and these orbits partition $\mathbb{Z}^{\omega} \setminus C$. Since there are $2^{\aleph_0}$-many of them, we enumerate them as $\{O_{\alpha}: \alpha < 2^{\aleph_0}\}$. Pick a representative $u_{\alpha} \in O_{\alpha}$ for every $\alpha < 2^{\aleph_0}$, so that $O_{\alpha} = \bigcup_{k \in \mathbb{Z}} [u_{\alpha} + \overline{k}]$. By choice of $C$, the sets $[u_{\alpha} + \overline{k}]$ and $[u_{\alpha} + \overline{l}]$ are disjoint for $k \neq l$. For each $\alpha$, let $A_{\alpha}$ denote the union of the even iterates of $[u_{\alpha}]$ under $g$ and let $B_{\alpha}$ denote the union of the odd iterates:
\begin{align*}
A_{\alpha} & = \bigcup_{k \in \mathbb{Z}} [u_{\alpha} + \overline{2k}] \\
B_{\alpha} & = \bigcup_{k \in \mathbb{Z}} [u_{\alpha} + \overline{2k+1}].
\end{align*}
Observe that $A_{\alpha}$ and $B_{\alpha}$ are disjoint and each closed under tail-equivalence, and $g[A_{\alpha}] = B_{\alpha}$ for every $\alpha$. Let $A = \bigcup_{\alpha} A_{\alpha}$ and $B = \bigcup_{\alpha}B_{\alpha}$. Then likewise $A$ and $B$ are disjoint, closed under tail-equivalence, and $g[A] = B$. Thus $\mathbb{Z}^{\omega} \setminus C = A \cup B$ is a partition of $\mathbb{Z}^{\omega} \setminus C$ into two isomorphic suborders, which are everywhere isomorphic by virtue of the maps $f_{I} \circ g \circ f_I^{-1}$. Since $\mathbb{Z}^{\omega} \setminus C \cong \mathbb{I}$, we are done. \qed
\,\ \\

We note that the same construction yields examples of homogeneous everywhere isomorphic sets $A, B$ of any infinite cardinality $\kappa \leq 2^{\aleph_0}$. Simply take $A = \bigcup_{\alpha < \kappa} A_{\alpha}$ and $B = \bigcup_{\alpha < \kappa} B_{\alpha}$. 

\begin{center} * \,\,\,\,\,\,\,\,\,\, * \,\,\,\,\,\,\,\,\,\, *
\end{center}

We end with some further discussion on homogeneity, everywhere isomorphism, and two-set partitions of $\mathbb{R}$. 

If $A \subseteq \mathbb{R}$ is homogeneous, then the order type of $A$ is unchanged by adding finitely many points to $A$. To see this, fix $x \in \mathbb{R}$. We claim $A \cup \{x\} \cong A$. If $x \in A$, there is nothing to show. If $x \not \in A$, then since $A$ is homogeneous and missing at least one point, it misses at least one point in every interval. That is, $\mathbb{R} \setminus A$ is dense in $\mathbb{R}$ (we say, $A$ is \emph{codense}). Since $A$ is homogeneous, it is also dense in $\mathbb{R}$. Thus we may pick a sequence $\ldots < x_{-1} < x_0 < x_1 < x_2 < \ldots$ of real numbers such that $x_0 = x$, $x_i \in A$ for all $i > 0$, and $x_i \in \mathbb{R} \setminus A$ for all $i < 0$. We assume that the sequence converges to $\infty$ on the right and $-\infty$ on the left since this is possible to arrange. Let $J_n$ denote the open interval $(x_n, x_{n+1})$ for all $n \in \mathbb{Z}$. We have $A \cap J_n \cong A$ for every $n$, and in particular $A \cap J_n \cong A \cap J_{n+1}$ for every $n$. Fix an isomorphism $f_n: A \cap J_n \rightarrow A \cap J_{n+1}$ for every $n$. Define a map $f: A \cup \{x\} \rightarrow A$ by the rule $f(y) = f_n(y)$ for all $y \in J_n$ and for all $n$, and $f(x_n) = x_{n+1}$ for all $n$. By construction, $f$ is an order-isomorphism. It follows that $A \cup X \cong A$ for any finite set of points $X \subseteq \mathbb{R}$. By a similar argument, one can show that if $A$ is homogeneous and $A \neq \mathbb{R}$, then $A \setminus X \cong A$ for any finite set of points $X$.

Is it true that for an arbitrary \emph{countable} set $X \subseteq \mathbb{R}$ we have $A \cup X \cong A$? This holds if $A$ is countable and homogeneous, since in this case both $A$ and $A \cup X$ are countable, dense, and without endpoints, and hence isomorphic to $\mathbb{Q}$. But it is not true in general for uncountable homogeneous $A$.

\theoremstyle{definition}
\newtheorem{cor3}[thm1]{Corollary}
\begin{cor3}\label{cor2} \,\ 
There is a homogeneous suborder $A \subseteq \mathbb{R}$ and countable set $X \subseteq \mathbb{R} \setminus A$ such that $A \cup X \not \cong A$. 
\end{cor3}

\begin{proof}
Let $\mathbb{R} \setminus \mathbb{Q} = A \cup B$ be our partition of the irrationals into everywhere isomorphic homogeneous sets. Along with this decomposition we have, for every open interval $I$, an isomorphism $f: \mathbb{R} \rightarrow I$ such that $f_I[A] = A \cap I$ and $f_I[B] = B \cap I$, as well as a global isomorphism $g: A \rightarrow B$. It follows $f_I[\mathbb{Q}] = \mathbb{Q} \cap I$ and so $f_I[A \cup \mathbb{Q}] = (A \cup \mathbb{Q}) \cap I$ as well, for every $I$. If there were an isomorphism $h: A \cup \mathbb{Q} \rightarrow A$, we would get the isomorphism $f_I \circ g \circ h \circ f_I^{-1}: (A \cup \mathbb{Q}) \cap I \rightarrow B \cap I$, for every $I$. But then $(A \cup \mathbb{Q})$ and $B$ would constitute a partition of $\mathbb{R}$ into everywhere isomorphic sets, a contradiction. Take $X = \mathbb{Q}$.
\end{proof}

Theorems 1 and 2 along with the corollary give us information about which subsets of $\mathbb{R}$ can be order-isomorphic to their complements. Suppose $A$ is a dense subset of $\mathbb{R}$ whose complement $B = \mathbb{R} \setminus A$ is also dense. Then if there is an isomorphism $f: A \rightarrow B$, it can be extended uniquely to a order-automorphism $f: \mathbb{R} \rightarrow \mathbb{R}$. Since this automorphism interchanges a set with its complement, it has no fixed points. If for a given $x$ we have $f(x) > x$, it follows that the iterate sequence $x < f(x) < f^2(x) < \ldots$ is unbounded to the right in $\mathbb{R}$ and the negative iterate sequence $\ldots < f^{-2}(x) < f^{-1}(x) < x$ is unbounded to the left, so that the same holds for \emph{every} $x \in \mathbb{R}$. The situation is symmetric if $f(x) < x$. We call such an automorphism of $\mathbb{R}$ \emph{irreducible}.

Let us suppose we have $f(x) > x$ for all $x$, since if not we may replace $f$ with its inverse, which still witnesses the isomorphism of $A$ and $B$. For any given $x$, the points in the iterate sequence $\ldots < f^{-1}(x) < x < f(x) < f^2(x) < \ldots$ alternate in their belonging to $A$ or $B$. If we consider, say, the interval $[0, f(0))$ and its partition into the two dense subsets $A \cap [0, f(0))$ and $B \cap [0, f(0))$, we have that $[f(0), f^2(0))$ is partitioned in the alternate way: $A \cap [f(0), f^2(0)) = f[B \cap [0, f(0))]$ and $B \cap [f(0), f^2(0)) = f[A \cap [0, f(0))]$. Continuing this line of thought, we see that the global partition $\mathbb{R} = A \cup B$ is determined by its restriction to the interval $[0, f(0))$ and the automorphism $f$. 

Thus the dense/codense subsets $A \subseteq \mathbb{R}$ that are isomorphic to their complements $B = \mathbb{R} \setminus A$ are in one-to-one correspondence with increasing irreducible automorphisms $f: \mathbb{R} \rightarrow \mathbb{R}$ equipped with a dense/codense partition of $[0, f(0))$. Since any automorphism of $\mathbb{R}$ is determined by its restriction to the rationals, a countable set, there are at most $2^{\aleph_0}$-many irreducible automorphisms. It is not hard to verify there are exactly $2^{\aleph_0}$-many. By inductively diagonalizing against all such automorphisms, it is possible to construct many examples of dense and codense suborders $A \subseteq \mathbb{R}$ that that are not isomorphic to their complements $B$.

However, for an arbitrary dense and codense $A \subseteq \mathbb{R}$, in general it does not seem easy to detect whether or not $A$ is isomorphic to $B = \mathbb{R} \setminus {A}$, at least without some record of $A$'s construction. But Theorem 1 tells us that it is easy if we know that $A$ is homogeneous. In this case, we always have that $A \not\cong B$. For if we had isomorphism between $A$ and $B$, we would have isomorphism everywhere by virtue of the homogeneity of $A$ and $B$, and hence a partition of $\mathbb{R}$ into two everywhere isomorphic sets. Saying this contrapositively, if $A$ \emph{is} isomorphic to its complement, then $A$ is not homogeneous: there must be some open interval $I$ for which $A \not\cong A \cap I$. 

On the other hand, Theorem 2 shows us that a homogeneous suborder of $\mathbb{R}$ can be within a countable set of being isomorphic to its complement, in two different senses. Let $\mathbb{I} = A \cup B$ be our decomposition of the irrationals into two homogeneous everywhere isomorphic sets, and let $X = A \cup \mathbb{Q}$ be the complement of $B$ in $\mathbb{R}$. Then $X$ is homogeneous and therefore necessarily not isomorphic to $B$, as we observed in the corollary.

If we wish to maintain homogeneity and get isomorphism, we can do this by deleting the countable set $\mathbb{Q}$ from $X$ to get $A$. But then $A \cup B$ is no longer a partition of $\mathbb{R}$. 

If we wish to get isomorphism and maintain a full partition of $\mathbb{R}$, this can also be achieved. Let $g: \mathbb{R} \rightarrow \mathbb{R}$ be our global isomorphism between $A$ and $B$ from the proof of Theorem 2, viewed as an order-automorphism of $\mathbb{R}$. It follows from the proof that $g$ is an irreducible automorphism. Partition $\mathbb{Q}$ into $g$-orbits, that is, sets of the form $\{\ldots, g^{-1}(x), x, g(x), g^2(x), \ldots\}$. From each orbit pick a point, and add its even iterates to $A$, and its odd iterates to $B$. Then by construction the resulting sets $A'$ and $B'$ are (globally) isomorphic, as witnessed by $g$, and both are within a countable set of being homogeneous. Moreover $A' \cup B'$ is a partition of $\mathbb{R}$. (By the discussion preceding the corollary, we could not have gotten away with a lesser modification, since any finite exchange of points from $X$ to $B$ would not have changed the order type of either set.) 

But in the process of passing from the near partition $A \cup B$ of $\mathbb{R}$ to an actual partition, we necessarily lose the homogeneity of our sets: there must now be an open interval $I$ such that $A' \cap I \not\cong A'$ and $B' \cap I \not\cong B'$.

\end{document}